\documentclass[11pt]{amsart}
\usepackage{amsmath,amssymb,latexsym,amsfonts,amscd}
\newtheorem{theorem}{Theorem}[section]
\newtheorem{proposition}[theorem]{Proposition}
\newtheorem{definition}[theorem]{Definition}

\newtheorem{lemma}[theorem]{Lemma}
\newtheorem{corollary}[theorem]{Corollary}

\newtheorem{question}[theorem]{Question}

\begin{document}

\title[invariants of local rings]{On
some local cohomology invariants of local rings}
\author{Gennady Lyubeznik}
\thanks{NSF support is gratefully acknowledged}
\address{Department of Mathematics, University of Minnesota, Minneapolis,
MN 55455}
\email{gennady@math.umn.edu}

\begin{abstract}
Let $A$ be a commutative Noetherian local ring containing a field of
characteristic $p>0$. The integer invariants $\lambda_{i,j}(A)$ have been
introduced in an old paper of ours. In this paper we completely describe
$\lambda_{d,d}(A)$ where $d={\rm dim}A$ in terms of the topology of
Spec$A$.
\end{abstract}

\maketitle
\section{Introduction}
All rings in this paper are commutative and Noetherian. Let $A$ be a local
ring that admits a surjection from a regular local
$n$-dimensional ring $R$ containing a field. Let
$I\subset R$ be the kernel of the surjection, let 
$\mathfrak m\subset R$ be the maximal ideal and let 
$k=R/\mathfrak m$ be the residue field of $R$. The Bass numbers 
$\lambda_{i,j}(A)={\rm
dim}_k{\rm Ext}^i_R(k, H^{n-j}_I(R))$ of the
local cohomology module $H_I^{n-j}(R)$ are all finite and depend only on
$A, i$ and
$j$, but  neither on $R$, nor on the surjection
$R\to A$. This has been proven in our  old paper \cite[Sec. 4]{L1} and these
invariants of local rings have since been studied by a number of authors
\cite{BB, GS, Ka1, Ka2, Wa}.

Let $d={\rm dim}A$. We have proven in our paper \cite[4.4i,ii,iii]{L1}
that
$\lambda_{i,j}=0$ if $j>d$ or $i>j$ while $\lambda_{d,d}(A)\ne 0$.
Kawasaki \cite{Ka1} and Walther \cite{Wa} have completely described
$\lambda_{d,d}(A)$ for $d\leq 2$, while Kawasaki \cite[Sec. 3, Cor.
1]{Ka2} proved that $\lambda_{d,d}(A)=1$ if $A$ is $S_2$. In
\cite[Sec. 7]{L3} we stated a question which we reproduce here in a more
precise form: 

\begin{question}
Is
$\lambda_{d,d}(A)$ equal to the number of connected components of the
graph $\Gamma_B$ where $B=\widehat{\hat A^{\rm sh}}$ is the completion
of the strict Henselization of the completion of $A$?
\end{question}

The graph $\Gamma_B$ for any local ring $B$ has been introduced by
Hochster and Huneke \cite[3.4]{HH}. We reproduce their definition:

\begin{definition} 
Let $B$ be a local ring. The graph $\Gamma_B$ is
defined as follows. Its vertices are the top-dimensional minimal primes
of $B$ (i.e. primes $P$ such that ${\rm dim}(B/P)={\rm dim}B$) and two
distinct vertices
$P$ and
$Q$ are joined by an edge if and only if the ideal $P+Q$ has height one.
\end{definition}

We expect that Question 1.1 has a positive answer and the goal
of this paper is to prove this in characteristic $p$.
Our main result is the following.

\begin{theorem}\label{maintheorem}
Let $A$ be a local $d$-dimensional ring containing a field of
characterisitc $p>0$. Let $B=\widehat{\hat A^{{\rm sh}}}$ be the
completion of the strict Henselization of the completion of $A$. Then
$\lambda_{d,d}(A)$ equals the number of connected components of the graph
$\Gamma_B$.
\end{theorem}

In particular, if $A$ is complete and has a separably closed residue field, then $B=A$ and so in this case we get the simpler statement that $\lambda_{d,d}(A)$ equals the number of connected components of the graph $\Gamma_A$.

Let $V$ be a
projective variety over a separably closed field and let $A$ be the local ring at the vertex of the
affine cone of some projective embedding of $V$. In \cite[p. 133]{L3} we
asked whether
$\lambda_{i,j}(A)$ depends only on $V$, $i$ and $j$ but not on the
embedding. Our Theorem 1.3 provides some
supporting evidence for a positive answer to this question by showing that
in characteristic $p>0$ the integer $\lambda_{d,d}(A)$ where $d={\rm
dim}A=1+{\rm dim}V$ indeed is independent of the emebedding and only
depends on the dimensions of the pairwise intersections of the irreducible
components of
$V$. 

Some of our arguments are characteristic-free and they might be
used in an eventual proof of the characteristic zero case. We have
collected them in Section 2. Theorem \ref{maintheorem} is proven in
Section 4.

\section{Characteristic-free results}
The main result of this section is Corollary 2.4 which reduces Question
1.1 to the case where $A$ is complete
with separably closed residue field (i.e. $A=B$) and $\Gamma_B$ is
connected, in which case we expect that $\lambda_{d,d}(A)=1$.

\begin{proposition}
Let $R$ be a regular local ring containing a field, let $I\subset R$ be an
ideal and let
$A=R/I$. Let $\Gamma_1,\Gamma_2,\dots,\Gamma_t$ be the connected
components of
$\Gamma_A$. Assume ${\rm height}I=h$. Let $I_i$ be the intersection
of the minimal primes of $I$ that are the vertices of $\Gamma_i$. Then
$H^h_I(R)\cong \oplus_iH^h_{I_i}(R)$.
\end{proposition}

\emph{Proof.} Let $I_0$ be the intersection of the minimal primes of $I$
of height bigger than $h$, so the radical of $I$ is $\cap_{i\geq
0}I_i$. Let $J=\cap_{i=0}^{i=t-1}I_i$. Then $I=J\cap I_t$ up to radical
and Mayer-Vietoris yields an exact sequence $$H^h_{J+I_t}(R)\to
H^h_J(R)\oplus H^h_{I_t}(R)\to H^h_I(R)\to H^{h+1}_{J+I_t}(R).$$ Since
every minimal prime of $I_0$ has height bigger than $h$ and does not
contain any minimal prime of $I_t$, the height of $I_0+I_t$ is bigger
than $h+1$. The height of $I_i+I_t$ for $1\leq i\leq t-1$ also is bigger
than $h+1$ because these ideals are intersections of primes corresponding
to the vertices of two different connected components of
$\Gamma_A$. Hence the ideal $J+I_t$ has height bigger than $h+1$ and
therefore $H^h_{J+I_t}(R)\cong H^{h+1}_{J+I_t}(R)\cong 0$. Thus the above
exact sequence implies an isomorphism $H^h_I(R)\cong H^h_J(R)\oplus
H^h_{I_t}(R)$. Now we are done by induction on $t$ considering that
$H^h_{I_0}(R)=0$ because the height of every minimal prime of $I_0$
is bigger than $h$.\qed

\medskip

\begin{lemma}\label{lemma}
Let $R$ be a regular local ring containing a field, let $I\subset R$ be an
ideal, let $A=R/I$, let $\mathfrak m\subset R$ be the maximal ideal of
$R$ and let $s=\lambda_{i,j}(A)$. Then $H^i_{\mathfrak
m}(H^{n-j}_I(R))\cong E^s$, a direct sum of $s$ copies of $E$ where $E$
is the injective hull of the residue field $k=R/\mathfrak m$ in the
category of
$R$-modules.
\end{lemma}

\emph{Proof.}
By
\cite[3.6a]{L1}, $H^i_{\mathfrak m}(H^{n-j}_I(R))$ is an injective
$R$-module. Since it is supported only on the maximal ideal $\mathfrak
m$, we have that  
$H^i_{\mathfrak m}(H^{n-j}_I(R))\cong E^s$, a direct sum of $s$ copies of
$E$, where $s$ is some finite or infinite cardinal. We quote \cite[1.4]{L1}:

{\it Let  $P$ be a prime of $R$ and let $M$ be an $R$-module such that
$(H^i_P(M))_P$ are injective for all $i$. Let $K(R/P)$ be the fraction
field of $R/P$. Then  
${\rm Ext}^i_{R_P}(K(R/P),M_P)\cong {\rm
Hom}_{R_P}(K(R/P),(H^i_P(M))_P)$.}

We set $M=H^{n-j}_I(R)$ and $P=\mathfrak m$. Since \hbox{$H^i_{\mathfrak
m}(H^{n-j}_I(R))\cong E^s$} while $K(R/P)=k$ and 
Hom$_R(k,E)\cong k$,
 we get
$\lambda_{i,j}(A)={\rm dim}_k{\rm
Hom}_R(k,E^s)=s$. \qed

\medskip

If a local ring $A$ containing a field does not admit a surjection from a
regular local ring, one sets
$\lambda_{i,j}(A)\stackrel{\rm def}{=}\lambda_{i,j}(\hat A)$ where $\hat
A$ is the completion of $A$ with respect to the maximal ideal. A
complete local ring containing a field always admits a surjection from a
complete regular local ring containing a field. In this way
$\lambda_{i,j}(A)$ are defined for every local ring
$A$ containing a field and one always has
$\lambda_{i,j}(A)=\lambda_{i,j}(\hat A)$
\cite[pp. 53-54]{L1}. For this reason we do not include the assumption
that $A$ is a surjective image of a regular local ring in Theorem
\ref{maintheorem} and in the following proposition.

\begin{proposition}
Let $A$ be a local ring containing a
field and let 
$B=\widehat{\hat A^{\rm sh}}$ be the completion of the strict
Henselization of the completion of
$A$. Then  $\lambda_{i,j}(A)=\lambda_{i,j}(B)$.
\end{proposition}

\emph{Proof.} Since $\lambda_{i,j}(A)=\lambda_{i,j}(\hat A)$, we can
assume that
$A$ is complete, i.e. there is a surjection $R\to A$ from a regular local
ring
$R$ containing a field. Let $I$ be the kernel of the surjection and let
$n={\rm dim}R$. By Lemma
\ref{lemma}, $H^i_{\mathfrak m}(H^{n-j}_I(R))\cong E^s$ where
$\lambda_{i,j}(A)=s$. 

Let
$T=\widehat{\hat R^{\rm sh}}$ be the completion of the strict
Henselization of the completion of $R$. Then $T/IT=B$. Since $T$ is flat
over $R$, we have that $H^i_{\mathfrak mT}(H^{n-j}_{IT}(T))\cong
T\otimes_RH^i_{\mathfrak m}(H^{n-j}_I(R))\cong (T\otimes_RE)^s$. It is
 enough to prove
that $T\otimes_RE\cong E_T$ where $E_T$ is the injective hull of the
residue field of $T$ in the category of $T$-modules for then
$H^i_{\mathfrak mT}(H^{n-j}_{IT}(T))\cong E_T^s$ and
$\lambda_{i,j}(B)=s$ by Lemma
\ref{lemma}. 

But $E\cong H^{n}_{\mathfrak m}(R)$ and $E_T\cong H^{n}_{\mathfrak
mT}(T)$ since $R$ and $T$ are regular, $\mathfrak mT$ is the maximal
ideal of $T$ and $n={\rm dim}R={\rm dim}T$. Considering that
$H^{n}_{\mathfrak mT}(T)=H^{n}_{\mathfrak
mT}(T\otimes_RR)=T\otimes_RH^{n}_{\mathfrak m}(R)$, we are
done.\qed

\smallskip

\begin{corollary}
Let $A$ be a local ring of dimension $d$ containing a field and let 
$B=\widehat{\hat A^{\rm sh}}$
be the completion of the strict Henselization of the completion of $A$. 
Let $\Gamma_1,\Gamma_2,\dots,\Gamma_r$ be the connected
components of
$\Gamma_B$. Let $I_j$ be the intersection
of the minimal primes of $B$ that are the vertices of $\Gamma_i$. Let
$B_j=B/I_j$. Then

(a)$\lambda_{i,d}(A)=\sum_{j=1}^{j=r}\lambda_{i,d}(B_j)$ for every $i$.

(b) Question 1.1 has a positive answer for $A$ (i.e.
$\lambda_{d,d}(A)=r$) if and only if it has a positive answer for every
$B_j$ (i.e. if and only if
$\lambda_{d,d}(B_j)=1$ for every $j$). 
\end{corollary}

\emph{Proof.} By Proposition 2.3, $\lambda_{i,j}(A)=\lambda_{i,j}(B)$.
Since $B$ is complete, there is a surjection $T\to B$ from a regular
local ring $T$ containing a field. Let $\tilde I$ be the kernel of the
surjection and let $\tilde I_j\subset T$ be the preimage of $I_j$, so
$T/\tilde I_j=B_j$. By Proposition 2.1, $H^{n-d}_{\tilde I}(T)\cong
\oplus_jH^{n-d}_{\tilde I_j}(T)$. Hence ${\rm Ext}^i_T(k,
H^{n-d}_{\tilde I}(T))\cong\oplus_j {\rm Ext}^i_T(k,
H^{n-d}_{\tilde I_j}(T))$ where $k$ is the residue field of $T$. Hence the
dimensions over
$k$ of the two sides of this equation are the same. This implies (a). And
(b) is immediate from (a).
\qed

\medskip

We note that each $B_j$ is complete, reduced, equidimensional, has a
separably closed residue field and $\Gamma_{B_j}$ is connected (since
$\Gamma_{B_j}=\Gamma_j$). Thus Corollary 2.4 reduces Question 1.1 to
rings of this type.

Finally, it is worth pointing out that the graph
$\Gamma_B$ where $B=\widehat{\hat A^{\rm sh}}$ is realized by a
substantially smaller ring than 
$\widehat{\hat A^{\rm sh}}$. Namely, let
$k\subset
\hat A$ be a coefficient field. It follows from \cite[4.2]{HL} that there
exists a finite separable field extension $K$ of $k$ such that the graphs
$\Gamma_{K\otimes_k\hat A}$ and $\Gamma_B$ are isomorphic.

\section{$F_T$-modules, $B\{f\}$-modules and the functor $\mathcal
H_{T,B}$.}\label{recall}
In this section we review 
some facts from our old paper \cite{L2} that are used in our
proof of Theorem \ref{maintheorem}.
Throughout this section $T$ is a complete local regular ring containing a
field of characteristic $p>0$.

\subsection{The Frobenius functor $F_T$.}Let
$T'$ be the additive group of
$T$ regarded as a
$T$-bimodule with the usual left $T$-action and with the right $T$-action
defined by 
$t't=t^pt'$ for all $t\in T, t'\in T'$. The Frobenius 
functor 
$$F_T=F:T\text{-mod}\to T\text{-mod}$$ 
of Peskine-Szpiro \cite[I.1.2]{PS} is defined by  
$$F(M)=T'\otimes_T M$$ 
$$F(M\stackrel{h}{\to}N)=(T'\otimes_RM\stackrel{\text{id}\otimes _T
h}{\longrightarrow}T'\otimes_TN)$$ for all $T$-modules $M$ and all
$T$-module homomorphisms $h$, where $F(M)$ acquires its $T$-module
structure via the left $T$-module  structure on $T'$. For a summary
of basic properties of the Frobenius functor
\cite[Remarks 1.0]{L2} may be consulted. 

\subsection{$F_T$-modules.}An $F$-module \cite[1.1]{L2} (more precisely,
an $F_T$-module) is a
$T$-module
$\mathcal M$ equipped with  a $T$-module 
isomorphism $\theta:\mathcal M\to F(\mathcal M)$ which we call the
structure  morphism of $\mathcal M$. 
A homomorphism of $F$-modules is a $T$-module homomorphism 
$f:\mathcal M\to \mathcal M'$ such that the 
following diagram commutes (where $\theta$ and $\theta'$ are the structure 
morphisms of $\mathcal M$ and $\mathcal M'$). 
$$
\CD
\mathcal{M} @>f>> \mathcal {M'}\\
@V\theta VV @VV\theta' V\\
F(\mathcal{M}) @>F(f)>> F(\mathcal{M}')
\endCD
$$

\subsection{Generating morphisms.}A generating morphism of an $F$-module
$\mathcal M$
\cite[1.9]{L2} is a
$T$-module  homomomorphism $\beta:M\to F(M)$, where $M$ is some
$T$-module, such that 
$\mathcal M$ is the limit of the inductive system in the top row of the 
commutative diagram  
$$
\CD
M @>\beta>> F(M) @>F(\beta)>> F^2 (M) @>F^2(\beta)>> \ldots\\
@V\beta VV @VV F(\beta) V @VV F^2(\beta) V\\
F(M) @>F(\beta)>> F^2 (M) @>F^2(\beta) >> F^3 (M) @>F^3(\beta) >> \ldots
\endCD
$$
and $\theta:\mathcal M\to F(\mathcal M)$, the structure isomorphism of
$\mathcal M$, is  induced by the vertical arrows in this 
diagram. Since the tensor product commutes with direct
limits, the limit of the inductive system of the bottom row  
is indeed $F(\mathcal M)$, so this definition of $\theta$ makes sense.

\subsection{Morphisms of $F$-modules in terms of generating
morphisms.}If $\beta:M\to F(M)$ and $\beta':M'\to F(M')$
are generating morphisms of $F$-modules $\mathcal M$ and $\mathcal
M'$ respectively, then any
$R$-module homomorphism \hbox{$h:M\to M'$} that makes the leftmost square
in the  diagram
$$
\CD
M @>\beta>> F(M) @>F(\beta)>> F^2 (M) @>F^2(\beta)>> \ldots\\
@Vh VV @VV F(h) V @VV F^2(h) V\\
M' @>\beta'>> F (M') @>F(\beta') >> F^2 (M') @>F^2(\beta') >> \ldots
\endCD
$$
commutative, makes the whole diagram commutative and the vertical arrows of 
this diagram induce an $F$-module homomorphism 
$\mathcal H:\mathcal M\to \mathcal M'$ \cite[1.10b]{L2}.

\subsection{Localization in terms of generating morphisms.}If $g\in T$ is
an element and
$\mathcal M$ is an
$F$-module, the localization $\mathcal M_g$ carries a natural structure of
$F$-module such that the natural localization map $\ell:\mathcal M\to
\mathcal M_g$ is a homomorphism of $F$-modules \cite[1.3b]{L2}. If
$\beta:M\to F(M)$ is a generating morphism of $\mathcal M$,
then $\beta\circ g^{p-1}:M\to F(M)$ is a generating morphism of
$\mathcal M_g$ 
\cite[1.10c]{L2} and the natural localization map $\ell:\mathcal M\to
\mathcal M_g$ is induced by the vertical arrows in the commutative diagram
$$
\CD
M @>\beta>> F(M) @>F(\beta)>> F^2 (M) @>F^2(\beta)>> \ldots\\
@V \text{mult by} VgV @V\text{mult by}V F(g)=g^p V @V\text{mult by}V
F^2(g)=g^{p^2} V\\ M @>\beta\circ g^{p-1}>> F (M) @>F(\beta\circ g^{p-1})
>> F^2 (M) @>F^2(\beta\circ g^{p-1}) >>
\ldots
\endCD
$$
Let $x\in F^s(M)$ (we mean the copy of $F^s(M)$ in the top row) and let
$\tilde x\in
\mathcal M$ be its image under the maps in the top row. Then the image of
$x\in F^s(M)$ (here we mean the copy of $F^s(M)$ in the bottom
row) in
$\mathcal M_g$ under the maps in the bottom row
is $\tilde x'=\frac{\ell(\tilde x)}{g^{p^s}}$. The $s$-th vertical map
which is the multiplication by $g^{p^s}$ on $F^s(M)$ induces a map
from the image of the top $F^s(M)$ in $\mathcal M$ 
to the image of the bottom $F^s(M)$ in $\mathcal M_g$ that sends
$\tilde x\in \mathcal M$ to $g^{p^s}\tilde x'=\frac{\ell(\tilde
x)}{1}\in \mathcal M_g$. 

\subsection{$B\{f\}$-modules.}Let $B$ be a surjective image of $T$ and
let 
$B\{f\}$ be the ring extension
of $B$ obtained by adjoining a variable $f$ subject to relations
$fb=b^pf$ for every $b\in B$. A (left) $B\{f\}$-module is a $B$-module
$N$ with a map
$f:N\to N$ called the action of the Frobenius on $N$ such that
$f(bx)=b^pf(x)$ for every
$x\in N$. We call $N$ a cofinite $B\{f\}$-module if it is cofinite (i.e.
Artinian) as a
$B$-module.

\subsection{The functor $\mathcal H_{T,B}$.}The surjective
homomorphism
$T\to B$ gives every 
$B$-module a structure of $T$-module and it gives every $T$-module 
annihilated by the kernel of the surjection a structure of $B$-module. If 
$N$ is a $B\{f\}$-module, we get a $T$-module homomorphism 
$$\gamma_N:F(N)=T'\otimes _TN\stackrel{t'\otimes x\mapsto
t'f(x)}{\longrightarrow} N.$$  Let $D={\rm Hom}_T(-,E)$ be the Matlis
duality functor where $E$ is the injective hull of the residue field of
$T$ in the category of
$T$-modules. Assume $N$ is a
cofinite $B\{f\}$-module. Applying
$D$ to the map
$\gamma_N$ and considering that there is a functorial $T$-module
isomorphism $\tau:D(F(N))\to F(D(N))$
\hbox{\cite[4.1]{L2}} we get a
$T$-module homomorphism 
$$\beta_N=\tau\circ D(\gamma_N):D(N)\to F(D(N)).$$ 
$\mathcal H_{T,B}(N)$ is the
$F$-module with generating morphism $\beta_N$. 
A homomorphism of $B\{f\}$-modules $\theta:N\to N'$ 
induces a homomorphism of $T$-modules \hbox{$D(\theta):D(N')\to D(N)$.} It
is straightforward to check that the left square in the diagram
$$
\CD
D(N')@>\beta_{N'}>>F(D(N'))@>F(\beta_{N'})>>F^2(D(N'))@>F^2(\beta_{N'})>>
\dots\\
@VV D(\theta)V  @VV F(D(\theta)) V  @VV F^2(D(\theta)) V\\
D(N)@>\beta_N>>F(D(N))@>F(\beta_N)>>F^2(D(N))@>F^2(\beta_N)>>\dots\\
\endCD
$$
is commutative, hence the whole diagram is commutative, hence 
 it induces an $F$-module homomorphism 
$\mathcal H_{T,B}(\theta):\mathcal H_{T,B}(N')\to \mathcal H_{T,B}(N)$
\cite[1.10b]{L2}. The functor $\mathcal H_{T,B}$ thus defined is a
functor from the category of cofinite $B\{f\}$-modules to the category of
$F$-modules supported on $V(I)$ where $I\subset T$ is the kernel of the
surjection $T\to B$. This functor is additive, contravariant and exact
\cite[4.2]{L2}.

\subsection{$H^i_{\mathfrak m}(B)$ as a $B\{f\}$-module.} The
local cohomology modules $H^i_{\mathfrak m}(B)$ have a natural
structure of $B\{f\}$-modules \cite[1.2a,b]{L2}. If the dimension
of $T$ is $n$, then $\mathcal H_{T,B}(H^i_{\mathfrak m}(B))\cong
H^{n-i}_I(T)$ \cite[4.8]{L2}.

\subsection{The stable part of a cofinite $B\{f\}$-module.}Since $B$ is
complete, it contains a coefficient field
$k\subset B$.  Let $N$ be a cofinite 
$B\{f\}$-module, let im$f^j$ be the set of elements of $N$ of the 
form $f^j(x)$ where $x\in N$, and let 
$k$-im$f^j$ be the $k$-vector subspace of $N$ spanned by im$f^j$. 
We set the stable part of $N$ to be the $k$-vector space 
$N_s=\cap_jk\text{-im}f^j$. This definition depends on the choice
of the coefficient field $k$; nevertheless we denote the stable
part simply by $N_s$. According to
\cite[4.9]{L2},
$N_s$ is finite-dimensional over $k$, $f:N_s\to N_s$ is injective and the
$k$-vector subspace of $N$ spanned by $f(N_s)$ coincides with $N_s$; this implies that if a set $\{y_1,\dots, y_q\}\subset N_s$ is a $k$-basis of $N_s$, then the set $\{f^v(y_1),\dots,f^v(y_q)\}\subset N_s$, for every $v$, is also a $k$-basis of $N_s$. The
dimension of $N_s$ as a $k$-vector space is independent of the choice of
the coefficient field $k$ \cite[4.11]{L2}.

If the dimension of the support of
$\mathcal H_{T,B}(N)$ is zero, it follows from \cite[4.10]{L2}
that
$\mathcal H_{T,B}(N)\cong E^r$ where $r={\rm dim}_kN_s$ 
(because in this case the maximal zero-dimensional
quotient of $\mathcal H_{T,B}(N)$ in the category of $F$-modules, denoted
$\mathcal L$ in
\cite[p. 110]{L2}, is
$\mathcal H_{T,B}(N)$ itself).

\section{Proof of Theorem \ref{maintheorem}}
We deduce Theorem \ref{maintheorem} from Corollary 2.4 by proving that
$\lambda_{d,d}(B_j)=1$ for every $B_j$. First we express
$\lambda_{i,j}(B)$ in terms of the Frobenius action on
$H^i_{\mathfrak m}(B)$ (Proposition 4.1 and Corollary 4.2). We then
use this to prove that if $B$ is one
the rings $B_j$ appearing in
Corollary 2.4, then $\lambda_{d,d}(B)=\lambda_{d,d}(S)$ where
$S$ is the $S_2$-ification of $B$ (Proposition 4.3) and conclude by
appealing to a result of Kawasaki to the effect that
$\lambda_{d,d}(S)=1$. 

Let $B$ be a complete local ring containing a field of characteristic
$p>0$. If
$N$ is a cofinite
$B\{f\}$-module and
$g\in B$, we let
$N(g)$ to be the cofinite $B\{f\}$-module defined as follows. The
underlying
$B$-module of $N(g)$ is $N$ and the action of the Frobenius on $N(g)$ is
defined by
$f(x)=g^{p-1}f({\rm id}(x))$ for every $x\in N(g)$ where 
id$:N(g)\to N$ is the identity map of the underlying
$B$-modules. Clearly $N(gg')=N(g)(g')$ for all $g,g'\in B$.

Since $B$ is complete,
there is a surjection $T\to B$ where $T$ is a complete regular local
ring containing a field. 

\begin{proposition}\label{mainlemma}
Let $B$ and $T\to B$ be as above. Let $\mathbf
g=\{g_1,\dots,g_d\}\subset B$ be a system of parameters of $B$ and let $N$
be a cofinite $B\{f\}$-module. Let
$K_{\bullet}(\mathbf g;N)$ be the Koszul complex of $N$ on $\mathbf g$,
namely
$$0\to K_d(\mathbf g;N)\to K_{d-1}(\mathbf g;N)\to\dots\to K_0(\mathbf
g;N)\to 0$$
where $K_r(\mathbf g;N)=\oplus_{1\leq i_1<i_2<\dots<i_r\leq
d}N(g_{i_1}\cdots g_{i_r})$ and the differentials are defined as follows.
The image of $x\in N(g_{i_1}\cdots g_{i_r})\subset K_r(\mathbf g;N)$
under the corresponding differential is $\sum_j(-1)^jg_{i_j}{\rm
id}_{i_1,\dots,\widehat{i_j},\dots,i_r}(x)\in K_{r-1}(\mathbf g;N)$ where 
$id_{i_1,\dots,\widehat{i_j},\dots,i_r}:N(g_{i_1}\cdots
g_{i_r})\to N(g_{i_1}\cdots\widehat{g_{i_j}}\cdots g_{i_r})$ is the
identity map on the underlying $B$-modules. 

(i) $K_{\bullet}(\mathbf g;N)$ is a complex in the category of
cofinite
$B\{f\}$-modules. Hence the cohomology modules $H_i(K_{\bullet}(\mathbf
g;N))$ of this complex are in a natural way cofinite $B\{f\}$-modules. 

(ii) $\mathcal H_{T,B}(H_i(K_{\bullet}(\mathbf g;N)))\cong H^i_{\mathfrak
m}(\mathcal H_{T,B}(N))$ where $\mathfrak m$ is the maximal ideal of $T$.

(iii) $H^i_{\mathfrak m}(\mathcal H_{T,B}(N))\cong E^r$ where
$r={\rm dim}_kH_i(K_{\bullet}(\mathbf g;N))_s$, $E$ is the injective
hull of $T/\mathfrak m$ in the category of $T$-modules and
$H_i(K_{\bullet}(\mathbf g;N))_s$ is the stable part of 
$H_i(K_{\bullet}(\mathbf g;N))$ with respect to some
coefficient field $k$ of $B$ (see (3.9)).
\end{proposition}

\emph{Proof.}
(i) The
map $\delta:N(g)\to N$ which is the multiplication by $g$ on the
underlying
$B$-module is a 
$B\{f\}$-module homomorphism because $\delta(f(x))=g\cdot g^{p-1}f({\rm
id}(x))=g^{p}f({\rm
id}(x))=f(g\cdot{\rm id}(x))=f(\delta(x))$ for every $x\in N(g)$.
Since $N(g_{i_1}\cdots g_{i_r})\cong
N(g_{i_1}\cdots\widehat{g_{i_j}}\cdots g_{i_r})(g_{i_j})$, the map 
$N(g_{i_1}\cdots g_{i_r})\to
N(g_{i_1}\cdots\widehat{g_{i_j}}\cdots g_{i_r})$ which is the
multiplication by $g_{i_j}$ is a morphism of $B\{f\}$-modules. Hence
the differentials of
$K_{\bullet}(\mathbf g;N)$ are morphisms of
$B\{f\}$-modules. This proves (i).

(ii) Let $\delta:N(g)\to N$ be the multiplication by $g$ as before
and let $\tilde g\in T$ be a lifting of $g$. The map $\delta$ may
be viewed as the multiplication by $\tilde g$. Associated to this
map is a commutative diagram like in (3.7):
$$
\CD
D(N)@>\beta_{N}>>F(D(N))@>F(\beta_{N})>>F^2(D(N))@>F^2(\beta_{N})>>
\dots\\
@VV D(\delta)V  @VV F(D(\delta)) V  @VV F^2(D(\delta)) V\\
D(N(g))@>\beta_{N(g)}>>F(D(N(g)))@>F(\beta_{N(g)})>>F^2(D(N(g)))@>
F^2(\beta_{N(g)})>>\dots\\
\endCD
$$
Since $\delta$ is the multiplication by $\tilde g$, we conclude that
$F^s(D(\delta))$ is the multiplication by $\tilde g^{p^s}$. Since the
underlying $B$-modules of $N(g)$ and $N$ are the same, $D(N(g))\cong
D(N)$. Identifying the underlying $B$-modules of $N(g)$ and $N$
and viewing the maps $\gamma_N$ and $\gamma_{N(g)}$
of (3.7) as two maps with the same source and target, we get
$\gamma_{N(g)}=\tilde g^{p-1}\gamma_N$, hence $\beta_{N(g)}=\beta_N\circ
\tilde g^{p-1}$. Putting all of this together we get that the above
commutative diagram takes the following form
$$
\CD
D(N)@>\beta_{N}>>F(D(N))@>F(\beta_{N})>>F^2(D(N))@>F^2(\beta_{N})>>
\dots\\
@V\text {mult by}V\tilde g V  @V\text {mult by}V \tilde g^p V  @V\text {
mult by}V \tilde g^{p^2} V\\
D(N)@>\beta_N\circ \tilde g^{p-1}>>F(D(N))@>F(\beta_N\circ
\tilde g^{p-1})>>F^2(D(N))@> F^2(\beta_N\circ \tilde g^{p-1})>>\dots\\
\endCD
$$
The limit of the inductive system in the top row is by definition
$\mathcal H_{T,B}(N)$. Comparing this commutative diagram with the one in
(3.5) we see that the limit of the inductive system of the bottom row is
$\mathcal H_{T,B}(N)_{\tilde g}$ and the vertical arrows induce the
natural localization map $\ell:\mathcal H_{T,B}(N)\to \mathcal
H_{T,B}(N)_{\tilde g}$. In other words, $\mathcal H_{T,B}(N(g))\cong
\mathcal H_{T,B}(N)_{\tilde g}$ and $\mathcal H_{T,B}(\delta)=\ell$.

Let
$\tilde {\mathbf g}=\{\tilde g_1,\dots,\tilde g_d\}\subset T$ be a
lifting of $\mathbf g$, i.e. $\tilde g_i$ is a lifting of $g_i$ for
every $i$. Considering that $N(g_{i_1}\cdots g_{i_r})\cong
N(g_{i_1}\cdots\widehat{g_{i_j}}\cdots g_{i_r})(g_{i_j})$ the above
implies that 
$\mathcal H_{T,B}$ transforms the map $N(g_{i_1}\cdots g_{i_r}) \to 
N(g_{i_1}\cdots\widehat{g_{i_j}}\cdots g_{i_r})$ which is the
multiplication by $g_{i_j}$ on the underlying $B$-module into the natural
localization map
$\mathcal H_{T,B}(N)_{\tilde g_{i_1}\cdots\widehat{\tilde g_{i_j}}\cdots
\tilde g_{i_r}}\to 
\mathcal H_{T,B}(N)_{\tilde g_{i_1}\cdots
\tilde g_{i_r}}$.

Hence the complex $\mathcal
H_{T,B}(K_{\bullet}(\mathbf g;N))$ is nothing but the \v Cech complex 
$C^{\bullet}(\tilde{\mathbf g};\mathcal H_{T,B}(N))$ of
$\mathcal H_{T,B}(N)$ with respect to $\tilde {\mathbf g}$, namely
$$0\to C^0(\tilde{\mathbf g};\mathcal H_{T,B}(N))\to 
C^1(\tilde{\mathbf g};\mathcal H_{T,B}(N))
\to\dots\to
C^d(\tilde{\mathbf g};\mathcal H_{T,B}(N))\to 0$$ 
where $C^r(\tilde{\mathbf
g};\mathcal H_{T,B}(N))=\mathcal H_{T,B}(K_r(\mathbf
g;N))=\oplus_{1\leq i_1<\dots<i_r\leq d}\mathcal
H_{T,B}(N)_{\tilde g_{i_1}\cdots \tilde g_{i_r}}$. The $i$-th cohomology
module of this
\v Cech complex is $H^i_J(\mathcal H_{T,B}(N))$ where $J$ is the ideal of
$T$ generated by $\tilde{\mathbf g}$ \cite[1.3]{L1}. Since the functor
$\mathcal H_{T,B}$ is exact, it commutes with the operation of taking the
cohomology of complexes. Hence $\mathcal H_{T,B}(H_i(K_{\bullet}(\mathbf
g;N)))\cong H^i(C^{\bullet}(\tilde {\mathbf g};\mathcal
H_{T,B}(N)))\cong H^i_J(\mathcal H_{T,B}(N))$. It remains to show that
$H^i_J(\mathcal H_{T,B}(N))\cong H^i_{\mathfrak m}(\mathcal H_{T,B}(N))$.

Since $\tilde g$ is a lifting of a system of parameters of $B$, the ideal
$J+I$ ($I$ is the kernel of the surjection $T\to B$) is
$\mathfrak m$-primary, i.e.
$H^i_{\mathfrak m}(M)\cong H^i_{I+J}(M)$ for every $T$-module $M$. The 
composition of functors $H^0_{\mathfrak m}(-)=H^0_J(H^0_I(-))$ yields a
spectral sequence $E_2^{p,q}=H^p_J(H^q_I(M))\Rightarrow
H^{p+q}_{\mathfrak m}(M)$. If
$M$ is supported on $V(I)$ then $H^q_I(M)=0$ for $q>0$ and $H^0_I(M)=M$,
i.e.
$E_2^{p,q}=0$ for $q>0$ and $E^{p,0}_2=H^p_J(M)$, so the spectral sequence
degenerates at $E_2$ and implies that $H^i_{\mathfrak m}(M)=H^i_J(M)$ for
all
$i$. Since
$\mathcal H_{T,B}(N)$ is supported on $V(I)$, this completes the proof of
(ii).

(iii) Since the dimension of the support of $H^i_{\mathfrak m}(\mathcal
H_{T,B}(N))$ is zero, we are done by (3.9) considering (ii). \qed

\medskip

The following corollary is crucial. It expresses $\lambda_{i,j}(B)$ in
terms of the Frobenius action on $H^j_{\mathfrak m}(B)$ without any
reference to the surjection $T\to B$.

\begin{corollary}\label{maincoro}
Let $B$ be a local ring containing a field of characteristic $p>0$. Let
$\mathbf g=\{g_1,\dots,g_d\}\subset B$ be a system of parameters of $B$
and let $\mathfrak m$ be the maximal ideal of $B$. Let
$K_{\bullet}(\mathbf g;H^j_{\mathfrak m}(B))$ be the Koszul complex of
$H^j_{\mathfrak m}(B)$ on $\mathbf g$. Let $H_i(K_{\bullet}(\mathbf
g;H^j_{\mathfrak m}(B)))$ be its $i$-th cohomology module. Then
$\lambda_{i,j}(B)={\rm dim}_kH_i(K_{\bullet}(\mathbf g;H^j_{\mathfrak
m}(B)))_s$.
\end{corollary}

\emph{Proof.} The local cohomology module $H^i_{\mathfrak m}(B)$ and
the action of the Frobenius on it remain the same after replacing $B$
by its completion with respect to $\mathfrak m$. Hence we may assume
that $B$ is complete with respect to its maximal ideal and therefore
admits a surjection $T\to B$ from a complete regular local ring $T$
containing a field of characteristic $p>0$. Let $I\subset T$ be the kernel
of this surjection and let $n={\rm dim}T$. Then $\mathcal
H_{T,B}(H^j_{\mathfrak m}(B))\cong H^{n-j}_I(T)$ by (3.8). Hence 
$H^i_{\mathfrak m}(H^{n-j}_I(T))\cong H^i_{\mathfrak m}(\mathcal
H_{T,B}(H^j_{\mathfrak m}(B)))\cong E^s$ where $s=\lambda_{i,j}(B)$ (by
Lemma 2.2 with
$A=B$ and
$R=T$). Now we are done by
Proposition 4.1iii with $N=H^j_{\mathfrak m}(B)$.\qed

\medskip 

If $B$ is complete, it has a canonical module. If in addition it is
reduced and equidimensional, the existence of a canonical module implies
that it has an
$S_2$-ification which we denote
$S$
\cite[2.7]{HH}. It is a ring extension of $B$ that is module-finite over
$B$ \cite[2.4]{HH}. A central result of \cite{HH} says that $S$ is a local
ring if and only if 
$\Gamma_B$ is connected \cite[3.6c,e]{HH}. 

\begin{proposition}
Let $B$ be a complete local ring containing a field of characteristic
$p>0$. Assume $B$ is reduced, equidimensional of dimension $d$, has a
separably closed residue field, and the graph $\Gamma_B$ is connected. Let
$S$ be the $S_2$-ification of $B$. Then
$\lambda_{i,d}(B)=\lambda_{i,d}(S)$ for every
$i$.
\end{proposition}

\emph{Proof.} Since $\Gamma_B$ is connected, $S$ is local, as is pointed
out above. Hence $\lambda_{i,d}(S)$ makes sense. 

Let $\mathfrak m_B$ and
$\mathfrak m_S$ be the maximal ideals of $B$ and $S$ respectively.
The ideal $\mathfrak m_BS$ of $S$ is $\mathfrak
m_S$-primary since $S$ is module-finite over $B$. This implies that 
$H^i_{\mathfrak m_B}(_BS)\stackrel{\text{by
\cite[5.7]{G}}}{\cong}_BH^i_{\mathfrak (m_B)S}(S)\cong _BH^i_{\mathfrak
m_S}(S)$ where the subscript $_B(-)$ means that the corresponding
$S$-module is viewed as a $B$-module via "restriction of scalars". 

Let
$Q=S/B$. The short exact sequence $0\to B\to S\to Q\to 0$ in the category
of $B$-modules yields an exact sequence $H^{d-1}_{\mathfrak m_B}(Q)\to
H^d_{\mathfrak m_B}(B)\to H^d_{\mathfrak m_B}(S)\to H^d_{\mathfrak
m_B}(Q)$. But for every $a\in S$ the ideal $\{b\in B|ba\in B\}$ of
$B$ has height at least two \cite[2.4]{HH}, so the natural inclusion $B\to
S$ becomes an isomorphism after localization at every prime of $B$ of
height at most one. Therefore the dimension of the $B$-module $Q$ is at
most
$d-2$. Hence $H^j_{\mathfrak m_B}(Q)=0$ for $j=d, d-1$ and the above exact
sequence implies an isomorphism of $B$-modules $H^d_{\mathfrak
m_B}(B)\cong H^d_{\mathfrak m_B}(S)\cong_BH^d_{\mathfrak m_S}(S)$ 
induced by the natural inclusion $B\to S$.

We claim that under this $B$-module isomorphism 
the natural action
of the Frobenius on $H^d_{\mathfrak m_B}(B)$ coincides with the natural
action of the Frobenius on $H^d_{\mathfrak m_S}(S)$, i.e. we have an
isomorphism of $B\{f\}$-modules, not just $B$-modules. Indeed, let
$\mathbf g=\{g_1,\dots,g_d\}\subset B$ be a system of parameters of $B$.
Then
$H^i_{\mathfrak m_B}(B)$ is the $i$-th cohomology of the \v Cech complex
$C^{\bullet}(\mathbf g;B)$ and the natural action of the Frobenius on
$H^i_{\mathfrak m_B}(B)$ is induced by the action of the Frobenius on
the complex $C^{\bullet}(\mathbf g;B)$, namely, if
$x\in B_{g_{i_1}\cdots g_{i_r}}\subset C^r(\mathbf g;B)$, then
$f(x)=x^p\in B_{g_{i_1}\cdots g_{i_r}}\subset C^r(\mathbf g;B)$. This
commutes with the
differentials and therefore induces an action of the Frobenius on
cohomology. Since
$S$ is module-finite over $B$, the set $\mathbf g$ is a system of
parameters for $S$ as well. Hence $H^i_{\mathfrak m_S}(S)$ is the $i$-th
cohomology of the \v Cech complex $C^{\bullet}(\mathbf g;S)$ and the
natural action of the Frobenius on $H^i_{\mathfrak m_S}(S)$ is induced by
the action of the Frobenius on the complex $C^{\bullet}(\mathbf
g;S)$, namely, if $x\in S_{g_{i_1}\cdots g_{i_r}}$, then $f(x)=x^p\in
S_{g_{i_1}\cdots g_{i_r}}$. The natural inclusion map
$B\to S$ induces a map of complexes
$C^{\bullet}(\mathbf g;B)\to _BC^{\bullet}(\mathbf g;S)$ which commutes
with the action of the Frobenius on both sides, as is easy to see. Hence
it induces a map on the $i$-th cohomology groups that commutes with the
action of the Frobenius. Since we have already seen that this map is an
isomorphism of $B$-modules for $i=d$, the claim is proven. 

Now let
$N=H_i(K_{\bullet}(\mathbf g;H^d_{\mathfrak m}(B)))$. By Corollary 4.2,
$\lambda_{i,d}(B)={\rm dim}_kN_{s,k}$ and
$\lambda_{i,d}(S)={\rm dim}_KN_{s,K}$ where $k\subset B$ and
$K\subset S$ are coefficient fields of
$B$ and
$S$ respectively and
$N_{s,k}$ and $N_{s,K}$ are the stable parts of $N$ with respect to $k$
and $K$ respectively, i.e. $N_{s,k}=\cap_jk\text{-im}f^j$ and
$N_{s,K}=\cap_jK\text{-im}f^j$. Thus it remains to prove that 
${\rm dim}_kN_{s,k}={\rm dim}_KN_{s,K}$. This equality trivially holds if
the residue fields of $B$ and $S$ coincide (which happens for example
if the reside field of $B$ is algebraically closed) for in this case
 a coefficient field of $B$ is
automatically a coefficient field of $S$, i.e. we may put $K=k$. In the
general case we are done by the following lemma. 
\qed

\begin{lemma}
Let $B$ and $S$ be as in Proposition 4.3. Let $N$ be a cofinite
$S\{f\}$-module. Let
$k\subset B$ and
$K\subset S$ be coefficient fields of
$B$ and
$S$ respectively. Let
$N_{s,k}$ and $N_{s,K}$ be the stable parts of $N$ with respect to $k$
and $K$ respectively, i.e. $N_{s,k}=\cap_jk\text{-im}f^j$ and
$N_{s,K}=\cap_jK\text{-im}f^j$. Then ${\rm dim}_kN_{s,k}={\rm
dim}_KN_{s,K}$.
\end{lemma}

\emph{Proof.} Viewing $B$ as a subring of $S$ we let $\tilde B=B+\mathfrak
m_S$ where
$\mathfrak m_S$ is the maximal ideal of $S$. Clearly $\tilde B$ is a
subring of $S$ containing $B$. Since $S$ is module-finite over $B$, so is
$\tilde B$. Hence $\tilde B$ is a complete local ring with maximal ideal
$\mathfrak m_{\tilde B}=\mathfrak m_S$ and $N$ is a cofinite
$\tilde B\{f\}$-module. Clearly,
$k\subset
\tilde B$ and $k$ is a coefficient field of $\tilde B$.

Let $\tilde k\subset K$ be the preimage of the image of $k$ in
$S/\mathfrak m_S$ under the natural map $k\cong B/\mathfrak m_B\to
S/\mathfrak m_S\cong K$. We claim $\tilde k\subset \tilde B$. Indeed,
for every $c\in \tilde k$ there is $c'\in k\subset B$ such that their
images in
$S/\mathfrak m_S$ are the same. Hence $c-c'\in \mathfrak m_S$. Therefore
$c=c'+(c-c')\in \tilde B$ which proves the claim.

Let $N_{s,\tilde k}$ be the stable part of $N$ with respect to $\tilde
k$. As is pointed out in (3.9), the dimension of the stable part is
independent of the choice of the coefficient field. Since
$k$ and
$\tilde k$ are two coefficient fields of the same ring $\tilde B$ and $N$
is a cofinite
$\tilde B\{f\}$-module, ${\rm dim}_kN_{s,k}={\rm dim}_{\tilde
k}N_{s,\tilde k}$. It remains to show that 
${\rm dim}_KN_{s,K}={\rm dim}_{\tilde
k}N_{s,\tilde k}$. 

Since $\tilde k\subset K$, we have that $\tilde k$-im$f^j\subset
K$-im$f^j$ for every
$j$, so $N_{s,\tilde k}\subset N_{s,K}$. Since $S$ is module-finite over
$B$, the residue field of $S$ is a finite field extension of the residue
field of $B$. Since the residue field of $B$ is separably closed, the
extension is purely inseparable. Thus $K\cong S/\mathfrak m_S$ is a finite
purely inseparable extension field of $\tilde k\cong k\cong B/\mathfrak
m$. Let $u$ be an integer such that $c^{p^u}\in \tilde
k$ for every $c\in K$.

Let dim$_{\tilde k}N_{s,\tilde k}=r$ and let $x_1,\dots, x_r\in
N_{s,\tilde k}$ be a
$\tilde k$-basis of
$N_{s,\tilde k}$. We claim $x_1,\dots, x_r$ are linearly independent
over $K$. Indeed, let $\sum_jc_jx_j=0$ be a linear dependency
relation where $c_j\in K$ are not all zero. Applying $f^u$ to this
relation we get $\sum_jc^{p^u}f^u(x_j)=0$ where $c_j^{p^u}\in \tilde
k$ for every $j$. Thus $f^u(x_1),\dots,f^u(x_r)$ are linearly dependent
over $\tilde k$. But this is impossible because the $\tilde k$-linear
span of these elements coincides with the $\tilde k$-linear span of
$f^u(N_{s,\tilde k})$ which according to (3.9) is just $N_{s,\tilde k}$
and $N_{s,\tilde k}$ has dimension $r$ over $\tilde k$. This proves the
claim and implies that $r={\rm dim}_{\tilde k}N_{s,\tilde k}\leq {\rm
dim}_KN_{s,K}$ since $N_{s,\tilde k}\subset N_{s,K}$. It remains to show
that ${\rm dim}_KN_{s,K}\leq r$.

Let $q={\rm dim}_KN_{s,K}$ and let $y_1,\dots, y_q\in N_{s,K}$ be a
$K$-basis of 
$N_{s,K}$. According to (3.9), $f^v(y_1),\dots,f^v(y_q)$ is a $K$-basis of $N_{s,K}$ for every $v$. Consequently $f^v(y_1),\dots,f^v(y_q)$ are linearly
independent over
$\tilde k$ which is a subfield of $K$. 

Let
$N'$ be the
$\tilde k$-linear span of $f^u(y_1),\dots,f^u(y_q)$. We claim that
$f$ sends $N'$ to itself. Indeed, $f(y_j)\in N_{s,K}$ for every $j$ since
according to (3.9), $f$ sends $N_{s,K}$ to itself. Hence
$f(y_j)=\sum_ic_iy_i$ where $c_i\in K$. Applying $f^u$ we get
$f^{u+1}(y_j)=\sum_ic_i^{p^u}f^u(y_i)\in N'$ since $c_i^{p^u}\in \tilde
k$ for every $i$. Since $f^{u+1}(y_1),\dots,f^{u+1}(y_q)$ span the
$\tilde k$-linear span of $f(N')$, the claim is proven. Since
$f^{u+1}(y_1),\dots,f^{u+1}(y_q)$ are linearly independent over
$\tilde k$, the $\tilde k$-linear span of $f(N')$ has dimension at
least $q={\rm dim}_{\tilde k}N'$, i.e. the $\tilde k$-linear span of
$f(N')$ coincides with $N'$. Hence $N'\subset \tilde k$-im$f^j$ for
every $j$, i.e. $N'\subset N_{s,\tilde k}$. This implies $q={\rm
dim}_{\tilde k}N'\leq {\rm dim}_{\tilde k}N_{s,\tilde k}=r$ and
completes the proof of the lemma. \qed

\medskip

Now let $B$ be one of the rings $B_j$ that appear in Corollary 2.4. Then
$B$ is complete, local, contains a field of characteristic $p>0$, has a
separably closed residue field, is reduced, equidimensional,
$d$-dimensional and the graph $\Gamma_B$ is connected. Hence it has an
$S_2$-ification which is a local ring $S$ and by Proposition 4.3,
$\lambda_{d,d}(B)=\lambda_{d,d}(S)$. Now we appeal to a result of
Kawasaki \cite[Sec 3, Prop.1]{Ka2} to the effect that if $S$ is a
$d$-dimensional
$S_2$ local ring that admits a surjection from a regular local ring
containing a field, then $\lambda_{d,d}(S)=1$ (see Proposition 4.5 below). This implies that
$\lambda_{d,d}(B_j)=1$ for every $j$, so Theorem
\ref{maintheorem} now follows from Corollary 2.4.\qed

\medskip

Finally, following a referee's suggestion, we include, for the reader's convenience, a proof of 
Kawasaki's result mentioned in the preceding paragraph.

\begin{proposition} \cite[Sec 3, Prop.1]{Ka2} Let $S$ 
be a $d$-dimensional
$S_2$ local ring that admits a surjection from a regular local ring $B$
containing a field. Then $\lambda_{d,d}(S)=1$.
\end{proposition}

\emph{Proof.} Let $B\to S$ be the surjection in question and let $I$ be its kernel.
Let $\mathfrak m$ be the maximal ideal of $B$.  
The composition of functors $\Gamma_{\mathfrak m}(-)=\Gamma_{\mathfrak m}(\Gamma_I(-))$
leads to the spectral sequence 
$$\mathcal E_2^{p,q}=H^p_{\mathfrak m}(H^q_I(B))\Longrightarrow H^{p+q}_{\mathfrak m}(B).$$
Let $n$ be the dimension of $B$. Since $B$ is regular, $H^n_{\mathfrak m}(B)\cong E$ 
where $E$ is the injective hull of the residue field of $B$ in the category of $B$-modules.
According to Lemma 2.2, all we need to show is that $H^d_{\mathfrak m}(H^{n-d}_I(B))\cong E$. But
$H^d_{\mathfrak m}(H^{n-d}_I(B))\cong \mathcal E_2^{d, n-d}$ and the abutment in total degree $n$ is $H^n_{\mathfrak m}(B)\cong E$. Hence 
it is enough to prove that all differentials going out of and coming into $\mathcal E^{d,n-d}_r$ are zero for $r\geq 2$, and all terms $\mathcal E^{p,q}_2$ with $p+q=n$ vanish unless $p=d$ and $q=n-d$. For this would imply that $H^d_{\mathfrak m}(H^{n-d}_I(B))$ is isomorphic to the abutment which is $E$.

The outgoing differenitials $d_r:\mathcal E_r^{d,n-d}\to \mathcal E_r^{d+r,n-d-r+1}$ are zero for all $r\geq 2$ since the target module $\mathcal E_r^{d+r,n-d-r+1}=0$. Indeed, this target module is a subquotient of $\mathcal E_2^{d+r,n-d-r+1}=H^{d+r}_{\mathfrak m}(H^{n-d-r+1}_I(B))$ which is zero since $H^{n-d-r+1}_I(B)$ is supported on $V(I)={\rm Spec}(B/I=A)$, hence the dimension of the support of $H^{n-d-r+1}_I(B)$ is at most $d={\rm dim}A$ while $d+r>d$ and $H^t_{\mathfrak m}(M)$, for any module $M$, vanishes provided $t$ is bigger than the dimension of the support of $M$.

The incoming differentials $d_r:\mathcal E_r^{d-r,n-d+r-1}\to \mathcal E_r^{d,n-d}$ are zero for all $r\geq 2$ since the source module $\mathcal E_r^{d-r,n-d+r-1}=0$. Indeed, this source module is a subquotient of $\mathcal E_2^{d-r,n-d+r-1}$, so it is enough to prove that the dimension of the support of $H^{n-d+r-1}_I(B)$, for all $r\geq 2$, is less than $d-r$, for this implies that $\mathcal E_2^{d-r,n-d+r-1}=H^{d-r}_{\mathfrak m}(H^{n-d+r-1}_I(B))=0$. 

Let $P\subset B$ be a prime ideal in the support of $H^{n-d+r-1}_I(B)$. Then $H^{n-d+r-1}_I(B)_P=H^{n-d+r-1}_{I_P}(B_P)\ne 0$. Since $B$ is a regular local ring, it is catenary, so if the dimension of $B/P$ is $\delta$, then the dimension of $B_P$ is $n-\delta$. If $\delta>d-r+1$, then the dimension of $B_P$ is $<n-d+r-1$, hence $H^{n-d+r-1}_{I_P}(B_P)=0$. Since $A=B/I$ is $S_2$ and catenary, it is equidimensional, hence the height of $IB_P$ is $n-d$. If $\delta=d-r+1$, then the dimension of $B_P$ is $n-d+r-1$ and $H^{n-d+r-1}_{I_P}(B_P)=0$ by the Hartshorne-Lichtenbaum local vanishing theorem \cite[3.1]{Ha} considering that the height of $I\widehat{B_P}$, where $\widehat{B_P}$ is the completion of the regular local ring $B_P$ with respect to its maximal ideal, is $n-d$, while the dimension of $\widehat{B_P}$ is $n-d+r-1>n-d$ (the latter inequality holds because $r\geq 2$). Finally, if $\delta=d-r$, then the dimension of $B_P$ is $n-d+r\geq n-d+2$, the height of $I_P$ is $n-d$, hence the dimension of $B_P/I_P$ is $\geq 2$. The fact that $A=B/I$ is $S_2$ impies that the depth of $B_P/I_P$ is $\geq 2$. Since $\widehat{\widehat{B_P/I_P}^{\rm{sh}}}$, the completion of the strict Henselization of the completion of $B_P/I_P$, is faithfully flat over $B_P/I_P$, the depth of $\widehat{\widehat{B_P/I_P}^{\rm{sh}}}$ also is $\geq 2$, which by \cite[2.1]{Ha1} implies that the punctured spectrum of $\widehat{\widehat{B_P/I_P}^{\rm{sh}}}$ is connected. These two facts, namely, dim$B_P/I_P\geq 2$ and the connectedness of the punctured spectrum of $\widehat{\widehat{B_P/I_P}^{\rm{sh}}}$, imply that $H^{n-d+r-1}_{I_P}(B_P)=0$ by \cite[Cor. 2.11]{O}, \cite[III, 5.5]{PS}, \cite[2.9]{HL} (the first two of these three references prove this result in characteristic 0 and $p>0$ respectively, while the third one gives a characteristic-free proof). All of this shows that if $\delta\geq d-r$, then $H^{n-d+r-1}_{I_P}(B_P)=0$, i.e. the dimension of the support of $H^{n-d+r-1}_{I}(B)$ is less than $d-r$. This completes the proof that all the incoming differentials are indeed zero.

It remains to show that if $p+q=n$, then $\mathcal E^{p,q}_2=0$ unless $p=d$ and $q=n-d$. Indeed, if $q<n-d$, then $H^q_I(B)=0$ since $B$ is regular and $q$ is smaller than the height of every minimal prime over $I$, which is $n-d$. Hence $\mathcal E^{p,q}_2=H^p_{\mathfrak m}(H^q_I(B))=0$ in this case. If $q>n-d$, i.e. $q=n-d+r-1$ for some $r\geq 2$, it has been shown in the preceding paragraph that the dimension of the support of $H^q_I(B)$ is less than $d-r=n-q-1=p-1$, hence $\mathcal E_2^{p,q}=H^p_{\mathfrak m}(H^q_I(B))=0$. \qed


\begin{thebibliography}{99}

\bibitem{BB}
M. Blickle and R. Bondu, {\it Local cohomology multiplicities in positive
characteristic}, preprint, (2004).

\bibitem{GS}
      R. Garcia Lopez and C. Sabbah,  
{\it Topological computation of local cohomology multiplicities},
Dedicated to the memory of Fernando Serrano. 
Collect. Math. {\bf 49} (1998), no. 2-3, 317--324.

\bibitem{G} A. Grothendieck,  Local cohomology, Lecture Notes in
Mathematics, v. {\bf 41}, Springer-Verlag, Heidelberg, 1967.

\bibitem{Ha1} R. Hartshorne, {\it Complete intersections and connectedness},
 Amer. J. Math. {\bf 84}  497--508  (1962).

\bibitem{Ha} R. Hartshorne, {\it Cohomological Dimension of Algebraic Varieties}, 
Ann. Math., {\bf 88} 403-450 (1968).

\bibitem{HH} M. Hochster and C. Huneke, {\it Indecomposable canonical
modules and connectedness},  Commutative algebra: syzygies,
multiplicities, and birational algebra (South Hadley, MA, 1992), 
197--208, Contemp. Math., {\bf 159}, Amer. Math. Soc., Providence, RI,
1994.

\bibitem{HL} C. Huneke and G. Lyubeznik, {\it On the vanishing of local
cohomology modules},  Invent. Math.  {\bf 102}  (1990),  no. 1, 73--93. 

\bibitem{Ka1}
        K.-I. Kawasaki,
                {\it On the Lyubeznik number of local cohomology modules},
 Bull. Nara Univ. Ed. Natur. Sci. {\bf 49} (2000) no. 2, 5-7.

\bibitem{Ka2}
       K.-I. Kawasaki,
               {\it On the highest Lyubeznik number}, Math. Proc. Cambr.
Phil. Soc.,  {\bf 132}  (2002),  no. 3, 409--417. 
 
\bibitem{L1} G. Lyubeznik, {\it Finiteness properties of local
cohomology modules (an application of $D$-modules to commutative
algebra)}  Invent. Math.  {\bf 113}  (1993),  no. 1, 41--55. 

\bibitem{L2} G. Lyubeznik, {\it $F$-modules: applications to local
cohomology and $D$-modules in characteristic $p>0$}, J. reine angew.
Math. {\bf 491} (1997), 65 - 130.

\bibitem{L3}  G. Lyubeznik, {\it A partial survey of local cohomology}, 
Local cohomology and its applications (Guanajuato, 1999),  121--154,
Lecture Notes in Pure and Appl. Math., {\bf 226}, Marcel Dekker, Inc., New
York, 2002.

\bibitem{O} A. Ogus, {\it Local Cohomological Dimension of Algebraic Varieties},
Ann. Math., {\bf 98}  327-365  (1973).

\bibitem{PS} C. Peskine and L. Szpiro, 
{\it Dimension projective finie et cohomologie locale. Applications ˆ la
dŽmonstration de conjectures de M. Auslander, H. Bass et A. Grothendieck},
 Inst. Hautes ƒtudes Sci. Publ. Math. No. {\bf 42} (1973), 47--119. 

\bibitem{Wa}
        U. Walther,
          {\it On the Lyubeznik numbers of a local ring}, 
Proceedings of the AMS, {\bf 129} (6) (2001) 1631-1634.

\end{thebibliography}
\end{document}